\documentclass[11pt]{article}
\usepackage{amsfonts,amsmath,amssymb,amsthm, amscd}
\usepackage[dvips]{epsfig,graphics}

\numberwithin{equation}{section}

\newtheorem{theorem}{Theorem}[section]
\newtheorem{prop}[theorem]{Proposition}
\newtheorem{lemma}[theorem]{Lemma}
\newtheorem{cor}[theorem]{Corollary}

\theoremstyle{definition}
\newtheorem{definition}[theorem]{Definition}
\newtheorem{example}[theorem]{Example}
\newtheorem{remark}[theorem]{Remark}

\newcommand{\D}{\Delta}

\def\<{{\langle}}
\def\>{{\rangle}}

\def\g{{\gamma}}
\def\Z{\mathbb Z}

\def\C{\mathbb C}

\def\l{{\lambda}}

\def\e{\epsilon}

\def\F{{\mathbb F}}

\begin{document}

\title{On Reciprocality of Twisted Alexander Invariants
 }

\author{
Jonathan A. Hillman\\ University of Sydney
 \and Daniel S. Silver \footnote{Supported in part by DMS-0706798.}
\\ University of South Alabama
\and
Susan G. Williams\footnote{Supported in part by DMS-0706798.}
\\ University of South Alabama
}

\maketitle

\begin{abstract} Given a knot and an ${\rm SL}_n \C$ representation of its group that is conjugate to its dual, the representation that replaces each matrix with its inverse-transpose, the associated twisted Reidemeister torsion is reciprocal. An example is given of a knot group and ${\rm SL}_3 \Z$ representation that is not conjugate to its dual for which the twisted 
Reidemeister torsion is not reciprocal. 
  \noindent \end{abstract}

\begin{footnote}{Mathematics Subject Classification:  
Primary 57M25.}\end {footnote}

\section{Introduction}  The Alexander polynomial $\D(t)$ of a knot $k$ can be computed from a diagram of $k$ or from a presentation of the knot group (see \cite{kaw}, for example). It is an integral Laurent polynomial, well defined up to multiplication by units $\pm t^i \in \Z[t^{\pm 1}]$, and it is usually normalized to be a polynomial with nonzero constant coefficient. 

It is well known that $\D(t)$ is {\it reciprocal} in the sense that
\begin{equation} \label{symmetry} \D(t^{-1}) \doteq  \D(t), \end{equation}
where $\doteq$ indicates equality up to multiplication by units. This is a consequence of
Poincar\'e duality of the knot exterior (see \cite{fox} for an alternative approach based on duality in the knot group).

In 1990 X.S.~Lin introduced a more sensitive invariant using information from nonabelian representations of the knot group \cite{lin}. Later, refinements were described by M. Wada \cite{wada} and others including P. Kirk and C. Livingston \cite{kl}, J. Cha \cite{cha}, and others. These twisted Alexander invariants have proven to be useful for a variety of questions about knots including questions about
concordance \cite{kl}, knot symmetry \cite{hln} and fibrations \cite{fvJAMS}. See \cite{friedl} for a survey. 

We briefly review the definition of perhaps the best-known twisted Alexander invariant. Let $k$ be a knot with exterior $X$, endowed with the structure of a CW complex.  We fix a Wirtinger presentation $\<x_0, x_1, \ldots, x_k \mid r_1, \ldots r_k\>$ for the knot group $\pi = \pi_1(X)$. 
Let $\phi: F_k \to \pi$ be the associated projection of the free group $F_k = \<x_0, x_1, \ldots, x_k \mid \>$ to $\pi$. It induces a ring homomorphism $\tilde \phi: \Z[F_k] \to \Z[\pi]$. 

Let $\e: \pi \to H_1(X; \Z) \cong \<t \mid \>$ be the abelianization mapping each $x_i$ to $t$. It induces a ring homomorphism $\tilde \e: \Z[\pi] \to \Z[t^{\pm 1}]$. 

Assume that $\g: \pi \to {\rm SL}_n \C$ is a linear representation.
Let $\tilde \g: \Z[\pi] \to M_n(\C)$ be the associated ring homomorphism to the algebra of $n\times n$ matrices over $\C$. We obtain a homomorphism 
\begin{equation} \tilde \g \otimes \tilde \e: \Z[\pi] \to M_n(\C[t^{\pm 1}]), \end{equation}
mapping $g$ to $t^{\e(g)} \g(g),$ that we denote more simply by $\Phi$. 

Let $M_{\g \otimes \e}$ denote the $k\times (k+1)$ matrix with $(i,j)$-component equal to the $n\times n$ matrix $\Phi(  {   {\partial r_i}\over{\partial x_j}  }) \in M_n(\C[t^{\pm 1}])$. Here ${\partial r_i}\over{\partial x_j}$ denotes Fox partial derivative.
Let $M^0_{\g \otimes \e}$ denote the $k\times k$ matrix obtained by deleting the column corresponding to $x_0$. We regard $M^0_{\g \otimes \e}$ as a $kn \times kn$ matrix with coefficients in $\C[t^{\pm 1}]$. 

\begin{definition} \label{wada} The {\it Wada invariant} $W_\g(t)$ is 
$${\det M^0_{\g \otimes \e} } \over {  \det \Phi(x_0 -1) }.$$ 
\end{definition}

When $\g$ is the trivial 1-dimensional representation, $M_{\g \otimes \e}^0$ is a matrix $M(t)$ that we call the {\it Alexander matrix}
of $k$. (This terminology is used, for example, in \cite{rolfsen}, but it is not standard.) The determinant of $M(t)$ is the (untwisted) Alexander polynomial $\D(t)$ of $k$. 

\begin{remark} \label{rational} The rational function $W_\g(t)$ need not be a polynomial. See \cite{wada}. \end{remark}

The matrix $M_{\g\otimes \e}$ represents a boundary homomorphism for a twisted 
chain complex 
\begin{equation} \label{ccplx} C_*(X; V[t^{\pm 1}]_\g) = (\C[t^{\pm 1}]\otimes_\C V)\otimes_\g C_*(\tilde X). \end{equation}
Here $V = \C^n$ is a vector space on which $\pi$ acts via $\g$, while $C_*(\tilde X)$ denotes the cellular chain complex of the universal cover $\tilde X$ with the structure of a CW complex that is lifted from $X$. The group ring $\Z[\pi]$ acts on the left via deck transformations.  On the other hand,  $\C[t^{\pm 1}]\otimes_\C V$ has the structure of of a right $\Z[\pi]$-module via 
$$(p \otimes v) \cdot g = (p t^{\e(g)})\otimes (v \g(g)),\ {\rm for\ } \g\in \pi.$$

\begin{remark} The homology group $H_1(X; V[t^{\pm 1}])$ of the chain complex (\ref{ccplx}) is a finitely generated $\C[t^{\pm 1}]$-module. Its 0th elementary divisor, $\D_\g(t)$, lately competes with $W_\g(t)$ for the  name ``twisted Alexander polynomial." In many cases they are equal; generally, $\D_\g(t)$ is $\det M_{\g \otimes \e}^0$ divided by a factor of $\det \Phi(x_0 - 1)$. See \cite{kl} or \cite{sw} for details. \end{remark} 

Let $\C(t)$ denote the field of rational functions. When $\det M_{\g \otimes \e}^0 \ne 0$, the chain complex

\begin{equation} \label{ccplx2} C_*(X; V(t)) = (\C(t)\otimes_\C V)\otimes_\g C_*(\tilde X) \end{equation}
is acyclic  \cite{kitano}, and hence the (Reidemeister) torsion $\tau_\g(t)$ is defined. In \cite{kl} it is shown that $\tau_\g(t)$ coincides with the Wada invariant $W_\g(t)$.

\begin{remark} Conjugating the representation $\g$ corresponds to a change of basis for $V$. It is well known that the invariants $W_\g(t), \D_\g(t)$ and $\tau_\g(t)$ are unchanged. \end{remark} 

T. Kitano used Poincar\'e duality to prove in \cite{kitano} that for orthogonal representations $\g: \pi \to {\rm SO}_n(\mathbb R)$, the torsion $\tau_\g(t)$ is reciprocal; that is, $\tau_\g(t^{-1}) \doteq \tau_\g(t)$. He asked whether reciprocality holds for general representations $\g: \pi \to {\rm SL}_n(\C)$.  The question appeared more recently in \cite{friedl}. 

Several years later, Kirk and Livingston showed in \cite{kl} that reciprocality holds whenever  $\g$ is unitary. In particular, it holds for all representations with finite image.

It is not difficult to find representations $\g: \pi \to {\rm GL}_n\C$ such
that $\tau_\g(t)$ is non-reciprocal. For example, consider the Wirtinger presentation $ \< x_0, x_1, x_2 \mid x_0x_1=x_2x_0, x_1x_2 = x_0x_1\>$ of the trefoil knot group $\pi$.
The assignment $x_i \mapsto X_i \in {\rm GL}_2\C$, such that
$$X_0=\begin{pmatrix} a & 0\\ 1 & 1 \end{pmatrix}, \quad  X_1=\begin{pmatrix} a & -(a^2-a+1)\\ 0 & 1 \end{pmatrix},\quad X_2=  X_1^{-1}X_0X_1$$
yields $\tau_\g(t) = at^2+1$. The question of reciprocality for representations in ${\rm SL}_n\C$ is more subtle. 

In Section \ref{examples} we show that reciprocality need not hold for general representations in ${\rm SL}_n\C$. The representations $\g$ that we consider have the property that the dual representation $\bar \g$, obtained by replacing each matrix $\g(g), g\in \pi$, by its inverse-transpose, is not conjugate to $\g$. 
We wish to thank Walter Neumann for suggesting to us that such a representation might yield non-reciprocal torsion.

In Section  \ref{condition} we prove that if a representation $\g: \pi \to {\rm SL}_n\C$ is conjugate to its dual, then the torsion $\tau_\g(t)$ is reciprocal. 

The authors wish to thank Kunio Murasugi for helpful suggestions.


\section{Examples} \label{examples} Any reciprocal even-degree integral polynomial $\D(t)$ such that $\D(1) = \pm 1$ arises as the Alexander polynomial of a knot (see \cite{kaw}, for example).  Let $f(t)$ be any monic integral polynomial with constant coefficient $-1$ and $f(1) = \pm 1$. Choose a knot $k$ with Alexander polynomial $\D(t)= f(t)f(t^{-1})$. 

Let $C$ be the companion matrix of $(t-1)f(t)$. Then $C \in {\rm SL}_n\Z$, where $\deg f = n-1$. 
Consider the cyclic representation $\g:\pi \to {\rm SL}_n \Z$ sending each generator $x_0, x_1, \ldots, x_k$ 
of a Wirtinger presentation of $\pi$ to $C$. We have 
\begin{equation} \label{torsion} \tau_\g(t) \doteq {{\det M^0_{\g \otimes \e} } \over {  \det \Phi(x_0 -1) } }\doteq {{\det M^0_{\g \otimes \e} } \over {  f(t^{-1})(t-1) }}.
\end{equation}

The matrix $M_{\g \otimes \e}^0$ can be obtained from the $(k\times k)$ Alexander matrix $M(t)$ by replacing each polynomial entry $\sum a_i t^i$ with the $(n\times n)$ block matrix $\sum a_i (tC)^i$. 
Since the $n\times n$ blocks commute, 
$$\det M_{\g \otimes \e}^0 = \prod_\l \det M(t \l),$$
where $\l$  ranges over the eigenvalues of $C$, that is, the roots of $(t-1)f(t)$ (see \cite{ksw} for details). Hence
$$\det M_{\g \otimes \e}^0 \doteq \prod_\l \D(t \l)\ =\ \D(t)\prod_{\l: f(\l)=0} f(t \l) f(t^{-1} \l^{-1}).$$
Since $\D(t)$ and $\det M_{\g \otimes \e}^0(t)$ are integral polynomials, so is 
$$g(t) = \prod_{\l: f(\l)=0} f(t \l) f(t^{-1} \l^{-1}).$$

\begin{lemma} \label{lemma} If $\deg f = 2$, then $g(t)$ is reciprocal. 
\end{lemma}

\begin{proof} Our assumptions about $f(t)$ imply that its roots have the form $\l, -\l^{-1}$, for some $\l \in \C\setminus \{0\} $. Then $g(t) = f(t \l) f(t^{-1} \l^{-1}) f(-t \l^{-1})f(-t^{-1}\l)$ while $g(t^{-1}) = f(t^{-1}\l)f(t \l^{-1})f(-t^{-1}\l^{-1})f(-t\l).$ Observe that $g(t)$ and $g(t^{-1})$ have the same roots: 
\begin{itemize} \item $f(t\l)$ and $f(-t^{-1}\l^{-1})$ have roots:  $t=1, -\l^{-2}$; \item $f(t^{-1} \l^{-1})$ and $f(-t\l)$ have roots: $t=-1, \l^{-2}$;
\item $f(-t\l^{-1})$ and $f(t^{-1}\l)$ have roots: $t=1, -\l^2$; 
\item $f(-t^{-1}\l)$ and $f(t \l^{-1})$ have roots: $t= -1, \l^2$. 
\end{itemize}
It follows  that $g(t^{-1}) = \alpha g(t)$, for some $\alpha \in \C\setminus \{0\}$. 
Letting $t=1$, we see that $\alpha =1$. Hence $g(t^{-1}) = g(t)$.

\end{proof} 

\begin{remark} The numerator   $\det M_{\g \otimes \e}^0$ of (\ref{wada}) is a polynomial invariant $D_\g(t)$ of $k$ (see \cite{sw}). Since $\D(t)$ is reciprocal, Lemma \ref{lemma} implies that $D_\g(t)$ is reciprocal whenever $\deg f =2$. Example \ref{ex2} below shows that this conclusion need not hold when $\deg f >2$. \end{remark}

\begin{prop} \label{propexample} Let $f(t)$ be a polynomial as above with degree $2$. If $f(t)$ is non-reciprocal, then $\tau_\g(t)$ is a non-reciprocal integral polynomial of the form $(t-1) h(t)$. \end{prop}

\begin{proof} From equation (\ref{torsion}),
\begin{equation} \label{torsion2} \tau_\g(t) \doteq {{   f(t)f(t^{-1})g(t)  } \over {f(t^{-1})(t-1)}}\doteq{ {f(t)g(t)} \over { t-1 }}.\end{equation}
Since $g(t)$ and $t-1$ are reciprocal but $f(t)$ is not, $\tau_\g(t)$ is non-reciprocal. To see that $\tau_\g(t)$ has the desired form, note that $(t-1)^2$ divides $g(t)$ since both factors $f(t \l), f(-t \l^{-1})$ of $g(t)$ vanish when $t=1$. 

\end{proof}

\begin{example} \label{ex1} Let $f(t) = t^2 - t-1$. Then 
$$C = \begin{pmatrix} 0 & 0 & -1\\ 1 & 0 & 0\\ 0 & 1 & 2 \end{pmatrix}.$$
Computation shows that $g(t) = (t-1)^2 (t+1)^2 (t^2 - 3 t +1)(t^2 + 3t + 1).$ 
By equation (\ref{torsion2}), 
$$\tau_\g(t) \doteq (t^2-t+1)(t-1)(t+1)^2(t^2-3t+1)(t^2+3t+1),$$
which is non-reciprocal. 
\end{example}

\begin{example} \label{ex2} Let $f(t) = t^3 - t-1$. Then 
$$C = \begin{pmatrix} 0 & 0 &0& -1\\ 1 & 0 & 0&0\\ 0 & 1 & 0&1\\ 0&0&1&1 \end{pmatrix}.$$
Computation shows that $g(t) = (t-1)^3(t^3-t-1)^2(t^3-t^2+2t-1)(t^6+3t^5+5t^4+5t^3+5t^2+3t+1).$ 
The polynomial $f(t)f(t^{-1})g(t)$ is the numerator $D_\g(t)$ of Wada's invariant (\ref{wada}). It is non-reciprocal. 

It is not difficult to see that for any cyclic representation, $D_\g(t) \doteq \D_\g(t)$ (see Section 3 of \cite{sw}) 
Hence this example shows that $\D_\g(t)$ can also be non-reciprocal.

\end{example}

\section{Sufficient condition for reciprocality}\label{condition}

If $\g: G \to {\rm GL}_n \F$ is a linear representation, then the {\it dual} (or {\it contragredient}) representation $\bar \g$ is defined by 
$$\bar \g(g) = {}^t\g(g)^{-1}, $$ where ${}^t$ denotes transpose. 

The following elementary lemma is included for the reader's convenience.

\begin{lemma} \label{dual} A representation $\g: G \to {\rm GL}_n \F$ is conjugate to its dual if and only if there exists a nondegenerate bilinear form $ (v, w) \mapsto \{v, w\}\in \F$ on $V$  such that $\{v\cdot g, w\cdot g\}=\{v, w\}$ for all $v, w \in V$ and $g \in G$. \end{lemma}

\begin{proof} Assume that $\bar \g$ is conjugate to $\g$. Then there exists a matrix
$A \in {\rm GL}_n \F$ such that $A^{-1} \g(g) A = {}^t\g(g)^{-1}$, for all $g \in G$. 
Define $\{v, w\} = v A\  {}^tw$. Since $A$ is invertible, the bilinear form is nondegenerate. It is easy to check that $\{v\cdot g, w\cdot g\}=\{v, w\}$ for all $v, w \in V$.

Conversely, assume that $\g$ preserves a nondegenerate bilinear form $(v, w) \mapsto \{v, w\}$. There exists an invertible matrix $A \in {\rm GL}_n \F$ such that $\{v, w\} = v A\ {}^tw.$ Since $\g$ preserves the form, we have $v \g(g) A\ {}^t\g(g)\ {}^tw =\{ v\cdot g, w\cdot g\} = \{v, w\} = v A\ {}^tw$, for all $v, w \in V, g \in G$. 
It follows that $\g(g)A\ {}^t \g(g) = A$ for all $g \in G$. Hence $A^{-1} \g(g) A= {}^t\g(g)^{-1}$, 
and so $\bar \g$ is conjugate to $\g$. 

\end{proof}

As before, let $k$ be a knot with group $\pi$. Assume that $\g: \pi \to {\rm SL}_n \F$ is a representation, where $\F$ is an arbitrary field. As above, $V = \F^n$ is a right $\Z[\pi]$-module via $v\cdot g = v \g(g),$ for all $v \in V$ and $\g \in \pi$. Let $W = \F^n$ with the dual $\Z[\pi]$-module structure given by $w\cdot g =w\  {}^t\g(t)^{-1}$. 

\begin{theorem} \label{thm} Assume that $\det M_{\g\otimes \e}^0 \ne 0.$ If $\g$ is conjugate to its dual representation $\bar \g$, then the torsion $\tau_\g(t)$ is reciprocal. \end{theorem}

\begin{proof}  The following argument is similar to those of \cite{kitano} and \cite{kl}. 

Recall that $X$ is the exterior of $k$, endowed with a CW cell structure.
Let $X'$ be the same space but with the dual cell structure. Let $\ \bar{}: \F(t) \to \F(t)$ be the involution induced by $t \mapsto t^{-1}$.

Assume that $\g: \pi \to {\rm SL}_n\F$ is a representation that is conjugate to its dual.
By Lemma \ref{dual} there exists a nondegenerate bilinear form $(v, w) \mapsto \{v\cdot g, w\cdot g\}$ such that $\{v\cdot g, w\cdot g\} = \{v, w\}$ for all $v, w \in V, g \in \pi$. Consider the twisted chain complexes
$$C_* = (\F(t)\otimes V) \otimes C_*(\tilde X),\ D_* = (\F(t)\otimes W)\otimes C_*(\tilde X', \partial \tilde X'),$$
where $\tilde X$ and $\tilde X'$ denote universal covering spaces of $X$ and $X'$, respectively. We abbreviate these by $V_{\g\otimes \e}\otimes C_*(\tilde X)$ and $V_{\bar \g\otimes  \e}\otimes C_*(\tilde X)$, respectively.

Define a bilinear pairing $C_q \times D_{3-q} \to F(t)$ by 
\begin{equation} \label{form} \<p \otimes v \otimes z_1, q \otimes w \otimes z_2\> = \sum_{g \in \pi} (z_1 \cdot g z_2) p \bar q \{v\cdot g, w\},\end{equation}
where $z_1 \cdot gz_2$ is the algebraic intersection number in $\Z$ of cells $z_1$ and $gz_2$. We extend linearly. 

The pairing induces a $\F(t)$-module isormorphism $D_{3-q} \to \overline{ {\rm Hom}}( C_q, \F(t))$, where $\overline{ {\rm Hom}}$ denotes the dual space with $(q\cdot h)(z) = \bar q(h(z))$, for all $q \in \F(t), z \in C_q$. Consequently, there exists a nondegenerate pairing
$H_q(X; V(t)) \times H_{3-q}(X', \partial X'; W(t)) \to \F(t)$. Since the torsion of $C_*$ is defined, by our hypothesis, the torsion of $D_*$ is too. 

Choose a basis $\{v_i\}$ over $\F$ for $V$ and lifts to $\tilde X$ of simplices of $X$. In this way, we obtain a preferred $\F(t)$-basis for $C_*$. Basis members have the form $1 \otimes v_i \otimes z_j$. We get a natural basis over $\F(t)$ for $D_*$ by picking a basis for $W$ that is dual to the basis for $V$ with respect to $\{, \}$, and choosing dual cells in $\tilde X'$ of the fixed lifts of simplices of $X$. As observed in \cite{kl}, the bases for $C_*$ and $D_*$ that we build are dual with respect to the bilinear form (\ref{form}). 

Let $\tau(X; V_{\g \otimes \e})$ denote the torsion of $C_*$. Similarly, let 
$\tau(X', \partial X'; V_{\bar \g \otimes \e})$ denote the torsion of $D_*$. Then $\tau(X; V_{\g\otimes \e}) =\tau(X', \partial X'; V_{\bar \g \otimes \bar \e})$ by Theorem $1'$ of \cite{m1}. Futhermore, 
\begin{align*}
\tau(X', \partial X'; V_{\bar \g \otimes \bar \e}) &=\tau(X, \partial X; V_{\bar \g \otimes \bar \e})\quad ({\rm by\ subdivision}) \\
	&= \tau(X, \partial X; V_{ \g \otimes \bar \e})\quad ({\rm since}\ \g\ {\rm is\ conjugate\ to}\ \bar \g) \\
	&=\bar \tau(X, \partial X; V_{\g \otimes \e})\\ 
	&=\bar \tau(X; V_{\g\otimes \e}).
\end{align*}
The last equality is a result of Lemma 2 of \cite{m2} and the fact that $\tau(\partial X; V_{\g \otimes \e})=1$ (see \cite{kl}).
 Hence 
 $$\tau_\g(t)= \tau(X; V_{\g\otimes \e}) = \bar \tau(X; V_{\g \otimes \e})= \bar \tau_\g(t).$$

\end{proof}
\begin{remark} \label{unitary} If $\F = {\mathbb R}$, and the bilinear form in Lemma \ref{dual} is positive-definite, then by considering a basis for $V$ that is  orthonormal with respect to the form, we see that $A$ is the identity matrix. In this case,  $\g(g) = {}^t\g(g)^{-1}$ for all $g \in G$, and hence $\g$ is conjugate to an orthogonal representation.
Similarly, if $\F = \C$ and the bilinear form is hermitian and positive-definite, $\g$ is conjugate to a unitary representation.

\end{remark}

\begin{cor} \label{symplectic} If $\g: \pi \to {\rm Sp}_{2n} \C$ is a symplectic representation,  then $\tau_\g(t)$ is reciprocal. \end{cor}

\begin{proof} The representation preserves the bilinear form given by $A=\begin{pmatrix} 0_n & I_n\\ -I_n & 0_n \end{pmatrix}$. \end{proof}

Since ${\rm Sp}_2\C = {\rm SL}_2 \C$, the following is immediate.

\begin{cor} \label{SL2}  If $\g$ is any representation of $\pi$ in ${\rm SL}_2\C$,  then $\tau_\g(t)$ is reciprocal. \end{cor}

Corollary \ref{SL2} shows that Example \ref{ex1} is, in a sense, the simplest possible.


 \bigskip

\noindent {\sl First author: School of Mathematics and Statistics F07, University of Sydney, NSW 2006 Australia} \\
{\sl E-mail:} jonh@maths.usyd.edu.au \medskip

\noindent {\sl Second and third authors:} Department of Mathematics and  Statistics, ILB 325, University of South Alabama, Mobile AL  36688 USA 
\\ \noindent {\sl E-mail:} silver@jaguar1.usouthal.edu; swilliam@jaguar1.usouthal.edu

\end{document}